\documentclass[12pt,leqno]{article}
\usepackage{amsmath, amsthm, amsfonts, amssymb, color}
\usepackage{mathrsfs}
\theoremstyle{definition}

\newcommand{\scr}[1]{\mathscr #1}
\definecolor{wco}{rgb}{0.5,0.2,0.3}

\numberwithin{equation}{section} \theoremstyle{remark}

\newcommand{\ua}{\uparrow}

\title{
{\bf  Harnack Inequality and Applications for Infinite-Dimensional GEM   Processes}
\footnote{Supported in part by NNSFC(11131003, 11431014), the 985 project, the Laboratory of Mathematical and  Complex Systems, and NSERC.}
}
\author{{\bf Shui Feng$^{b}$ and Feng-Yu Wang$^{a,c}$}\\
\footnotesize{$^a$School of Mathematical Sciences, Beijing Normal University, Beijing 100875, China}\\
 \footnotesize{$^b$Department of Mathematics and Statistics, McMaster University, Hamilton, L8S 4K1, Canada}\\
    \footnotesize{$^c$Department of Mathematics, Swansea University, Singleton Park, SA2 8PP, UK}\\
\footnotesize{  \tttext{wangfy@bnu.edu.cn}; \tttext{F.Y.Wang@swansea.ac.uk}; \tttext{shuifeng@univmail.cis.mcmaster.ca}} }

\begin{document}
\def\tttext#1{{\normalfont\ttfamily#1}}
\def\R{\mathbb R}  \def\ff{\frac} \def\ss{\sqrt} \def\B{\mathbf B}
\def\N{\mathbb N} \def\kk{\kappa} \def\m{{\bf m}}
\def\dd{\delta} \def\DD{\Delta} \def\vv{\varepsilon} \def\rr{\rho}
\def\<{\langle} \def\>{\rangle} \def\GG{\Gamma} \def\gg{\gamma}
  \def\nn{\nabla} \def\pp{\partial} \def\EE{\scr E}
\def\d{\text{\rm{d}}} \def\bb{\beta} \def\aa{\alpha} \def\D{\scr D}
  \def\si{\sigma} \def\ess{\text{\rm{ess}}}
\def\beg{\begin} \def\beq{\begin{equation}}  \def\F{\scr F}
\def\Ric{\text{\rm{Ric}}} \def\Hess{\text{\rm{Hess}}}
\def\e{\text{\rm{e}}} \def\ua{\underline a} \def\OO{\Omega}  \def\oo{\omega}
 \def\tt{\tilde} \def\Ric{\text{\rm{Ric}}}
\def\cut{\text{\rm{cut}}} \def\P{\mathbb P}
\def\C{\scr C}     \def\E{\mathbb E}\def\y{{\bf y}}
\def\Z{\mathbb Z} \def\II{\mathbb I}
  \def\Q{\mathbb Q}  \def\LL{\Lambda}\def\L{\scr L}
  \def\B{\scr B}    \def\ll{\lambda} \def\a{{\bf a}} \def\b{{\bf b}}
\def\vp{\varphi}\def\H{\mathbb H}\def\ee{\mathbf e}\def\x{{\bf x}}

\maketitle
\begin{abstract}  The dimension-free Harnack inequality and uniform heat kernel upper/lower bounds are derived for a class of  infinite-dimensional GEM processes, which was introduced in \cite{FW} to simulate the two-parameter GEM distributions.  In particular, the associated  Dirichlet form satisfies the super log-Sobolev inequality which strengthens the log-Sobolev inequality derived in \cite{FW}. To prove the main results, explicit Harnack inequality and super Poincar\'e inequality are established for the one-dimensional Wright-Fisher diffusion processes. The main tool of the study is the coupling by change of measures. \end{abstract} \noindent

   \vskip 2cm

\section{Introduction}

The GEM distribution appears in population genetics describing the distribution of age-ordered allelic frequencies (\cite{Ewens04}).   Due to the many computational friendly properties of the stick-breaking structure, the GEM distribution and various generalizations are widely used as prior distributions in Bayesian statistics (\cite{IshJames01}). Below we briefly recall a standard construction of the GEM random variables.

Let $\{U_n\}_{n\ge 1}$   be a sequence of independent beta random variables with corresponding parameters  $a_n>0$ and $b_n>0$, $n\ge 1$. Set
\begin{equation}\label{stick-breaking}
V_1=U_1, \ V_n=U_n \prod_{i=1}^{n-1}(1-U_i),\ \  n\geq 2.
\end{equation}

For any $n\geq 1$, the joint distribution of $(V_1, \ldots, V_n)$ is the generalized Dirichlet distribution defined in \cite{ConMoi69}. The law of ${\bf V}=(V_1,V_2,\ldots)$ is  a probability on the space $$\bar\DD_\infty:=\Big\{\x=(x_i)_{i\in\N}\in [0,1]^\N: \sum_{n=1}^\infty x_n\le 1\Big\}$$  equipped with the usual $\si$-field induced by the projections $\{\x\mapsto x_i: i\in \N\}.$ Let
\[
W_n =V_1+\ldots+V_n,\ \ \ n\ge 1.
\]
Then $W_n$ is monotonically increasing bounded above by $1$. If the parameters satisfy
\begin{equation}\label{dustless}
\sum_{i=1}^{\infty}\frac{a_i}{a_i +b_i}=\infty,
\end{equation}
then $1-W_n=(1-U_1)\cdots(1-U_n)$ converges monotonically to $0$ and the law of ${\bf V}$ becomes a probability on space
$$
\DD_\infty:=
\Big\{\x\in [0,1]^\N: \sum_{i\ge 1}x_i=1\Big\}.$$

If $a_i=1-\alpha,b_i=\theta+i\alpha$ for a pair of parameters $0\leq \alpha<1, \theta+\alpha>0$, then the law of ${\bf V}$ is the well known two-parameter GEM distribution.
The GEM distribution with parameter $\theta$, coined by Ewens and named after Griffiths, Engen, and McCloskey, corresponds to $\alpha=0$.  Under assumption \eqref{dustless}, the representation \eqref{stick-breaking} is also known as the stick-breaking model.

To simulate the GEM distributions using Markov processes, a class of infinite-deimensional diffusion processes on $\bar\DD_\infty$  have been constructed in \cite{FW}. It was proved in \cite{FW} that these processes are symmetric with respect to  the corresponding GEM distributions and satisfy the  log-Sobolev inequality, so that they converge to the GEM distributions exponentially in both entropy and $L^2$. In this paper we derive some    stronger properties on these processes, including the uniform  heat kernel upper/lower bounds  and super log-Sobolev inequalities.  The main idea of the study goes back to \cite{W97} using the dimension-free Harnack inequality, and the main tool to establish the Harnack inequality  is the coupling by change of measures developed from \cite{ATW}, see the recent monograph \cite{Wbook2} for a brief theory on coupling by change of measures and applications.

To recall the GEM processes constructed in \cite{FW}, let $\{a_i, b_i\}_{i\ge 1}$ be strictly positive numbers. Then the corresponding GEM process  is generated by the following second-order differentiable operator on $\bar \DD_\infty$ (note that the factor $\ff 1 2$ in the diffusion term is missed in \cite{FW}):
$$\L ({\bf x})= \ff 1 2  \sum_{i,j= 1}^\infty A_{ij}({\bf x})\ff{\pp^2}{\pp x_i\pp x_j} +
\sum_{i=1}^\infty C_i({\bf x})\ff{\pp}{\pp x_i},\ \ \x=(x_1,x_2,\cdots)\in \bar\DD_\infty, $$ where
\beg{equation*}\beg{split} & A_{ij}({\bf x}):= x_i x_j
\sum_{k=1}^{i\land j} \ff{(\dd_{ki}(1-\sum_{l=1}^{k-1}x_l)-x_k )
(\dd_{kj}(1-\sum_{l=1}^{k-1}x_l)-x_k )}{x_k(1-\sum_{l=1}^k
x_l)},\\
&C_i({\bf x}):= x_i\sum_{k=1}^i
\ff{(\dd_{ik}\big(1-\sum_{l=1}^{k-1}x_l\big)-x_k ) (a_k
\big(1-\sum_{l=1}^{k-1}x_l\big)- (a_k+b_k)x_k )}{x_k(1-\sum_{l=1}^k
x_l)}.\end{split}\end{equation*} Here and in what follows, we set $
\sum_{i=1}^0 =0$ and $\prod_{i=1}^0 =1$ by conventions. Obviously, $A_{ij}(\x)$ and $C_i(\x)$ are well defined if $\x=(x_i)_{i\in\N}$ satisfies $\sum_{i=1}^n x_i<1$ for all $n\in\N$. By setting $\ff 0 0=1$, they are defined on the whole space $\bar\DD_\infty$.

The  diffusion process generated by $\L$ on $\bar\DD_\infty$ is constructed in \cite{FW} using   the one-dimensional Wright-Fisher diffusion processes, which solve the following SDEs on $[0,1]$   for $i\ge 1$:
\beq\label{1.1} \d x_i(t)= \{a_i-(a_i+b_i)x_i(t)\}\d t +\ss{x_i(t)(1-x_i(t))}\,\d B_i(t), \end{equation}
where $\{B_i(t)\}_{i\ge 1}$ are independent one-dimensional Brownian motions. By \cite[Theorem 3.2]{IW} with $\si(x):=\ss{x(1-x)}1_{[0,1](x)}$ and $b(x):= a_i-(a_i+b_i)x$, the equation has a unique strong solution which is a diffusion process on $[0,1].$
For any $\x=(x_i)_{i\in\N}\in [0,1]^\N$, let $X^\x(t)=(x_1(t),x_2(t),\cdots)$, where   $x_i(t)$ solves \eqref{1.1} with $x_i(0)= x_i \in [0,1].$
Let $\bar P_t^{\a,\b}$ be the corresponding Markov semigroup, i.e.
$$\bar P_t^{\a,\b} f(\x)= \E f(X^\x(t)),\ \ \ t\ge 0, f\in \B_b([0,1]^\N), \x\in [0,1]^\N,$$ where $\B_b(\cdot)$ denotes the set of all bounded measurable functions on a measurable space.

It is easy to see that $x_i(t)$ is reversible with respect to the beta distribution
$$\pi_{a_i,b_i}(\d x):= \ff{\GG(2a_i+2b_i)}{\GG(2a_i)\GG(2b_i)} 1_{[0,1]}(x) (1-x)^{2a_i-1}(1-x)^{2b_i-1}\d x.$$
Define the map $\Phi: [0,1]^\N\to \bar\DD_\infty$ by
$$\Phi(\x)= (\phi_1(\x),\phi_2(\x),\cdots),\  \   \phi_n(\x):= x_n \prod_{i=1}^{n-1}(1-x_i), \ \ n\ge 1, \x=(x_1,x_2,\cdots)\in [0,1]^\N.$$
Let $\Xi_{\a,\b}= \pi_{\a,\b}\circ \Phi^{-1},$ where $\pi_{\a,\b}:= \prod_{i\ge 1} \pi_{a_i,b_i}.$  It is clear that $\Xi_{\a,\b}$ includes the GEM distributions as special examples:  the one-parameter GEM distribution $\pi_\theta^{GEM} =\Xi_{\a,\b}$ for $a_i=\ff 1 2$ and $b_i= \ff\theta 2$, and the two-parameter GEM distribution $\pi_{\aa,\theta}^{GEM}= \Xi_{\a,\b}$ for $a_i=\ff{1-\aa}2$ and $b_i= \ff{\theta+\aa i} 2.$

To construct the GEM diffusion process using the map $\Phi$ and $x_i(t), i\ge 1,$  we observe that  $$\Phi: [0,1)^\N\to \tt\DD_\infty=\Big\{\x\in [0,1]^\N:\ \sum_{i=1}^n x_i<1\ \forall\ n\in\N\Big\}$$ is a bijection with inverse
$$\Psi(\x) =(\psi_1(\x),\psi_2(\x),\cdots),\ \ \psi_n(\x):= \ff{x_n}{1-\sum_{i=1}^{n-1}x_i}\in [0,1),\ \ n\ge 1, \x\in \tt\DD_\infty.$$
Due to this fact, $\inf_{i\ge 1} b_i\ge \ff 1 2$ has been assumed in \cite{FW} so that    $x_i(t)\in [0,1)$ for all $t>0$ and $i\ge 1.$
In this case, for any $\x\in \Phi([0,1)^\N)$,  $Y^\x(t):= \Phi(X^{\Psi(\x)}(t))$ is a Markov process on $\tt\DD_\infty$.
 Moreover, according to   \cite[\S 3]{FW}, this Markov   process is generated by $\L$ on $\bar\DD_\infty$; that is,  the Markov semigroup
   \beq\label{F}P_t^{\a,\b}f(\x):= \E f(\Phi(X^{\Psi(\x)}(t)))=\bar P_t^{\a,\b}(f\circ\Phi)(\Psi(\x)),\ \ f\in \B_b(\tt\DD_\infty), t\ge 0, \x\in \tt\DD_\infty,\end{equation}
   where  $\B_b(\cdot)$ denotes the set of all bounded measurable real functions on a measurable space,
 is associated to the symmetric Dirichlet form $(\EE_{\a,\b},\D(\EE_{\a,\b}))$,
which is the closure of the following pre-Drichlet form on $L^2(\Xi_{\a,\b})$:
$$\EE_{\a,\b}(f,g):= -\int_{\bar\DD_\infty} f\L g\,\d\Xi_{\a,\b}= \ff 1 2 \int_{\bar\DD_\infty} \Big(\sum_{i,j\ge 1} a_{ij} (\pp_if)(\pp_jg)\Big)\,\d \Xi_{\a,\b}\ \ f,g\in \F C_b^\infty,$$
where $\F C_b^\infty$ is the set of all $C_b^\infty$-cylindrical functions on $[0,1]^\N.$

\

To extend the above construction for all $b_i>0$ for which $x_i(t)$ may hit $1$, we extend $\Psi$ to $\bar\DD_\infty$ by setting $\ff 0 0=1,$ i.e.
$\psi_n(\x)=1$ provided $\sum_{i=1}^{n-1} x_i=1$ (this implies $x_n=0$ for $\x\in \bar\DD_\infty$). Then
$$\Psi(\bar\DD_\infty)=E:=\big\{\x=(x_i)_{i\in\N}\in [0,1]^\N:\ \text{if}\ x_i=1\ \text{for\  some}\ i\in\N, \ \text{then}\ x_j=1\ \text{for\ all}\ j\ge i\big\},$$
and $\Phi: E\to \bar\DD_\infty$ is a bijection with inverse $\Psi$. In this case we  can prove that $P_t^{\a,\b}$ given in \eqref{F} for $\bar\DD_\infty$ in place of $\tt \DD_\infty$, i.e.
\beq\label{F2}P_t^{\a,\b}f(\x):= \E f(\Phi(X^{\Psi(\x)}(t)))=\bar P_t^{\a,\b}(f\circ\Phi)(\Psi(\x)),\ \ f\in \B_b(\bar\DD_\infty), t\ge 0, \x\in \bar\DD_\infty,\end{equation}
is also a Markov semigroup.
Indeed, since $\Phi\circ\Psi(\x)=\x$ for $\x\in \bar\DD_\infty,$   $P_0^{\a,\b}$ is the identity operator. Moreover, for any $t>0$ and any $\x\in [0,1]^\N$, we have
$$\P(X^\x(t)\notin E) \le \P(x_i(t)=1\ \text{for\ some}\ i\in \N)=0,$$ so that $\Psi\circ\Phi(X^\x(t))= X^\x(t)\ \P$-a.s. Thus,
by \eqref{F2} and the semigroup property of $\bar P_t^{\a,\b}$,
\beg{equation*}\beg{split}& P_t^{\a,\b} P_s^{\a,\b} f(\x) = \E\Big[(P_s^{\a,\b} f)\circ\Phi \big(X^{\Psi(\x)}(t)\big)\Big]\\
&= \E\Big[\big(\bar P_s^{\a,\b} \{f\circ\Phi\}\big)\circ\Psi\circ \Phi  \big(X^{\Psi(\x)}(t)\big)\Big]
 = \E\Big[ \bar P_s^{\a,\b} \{f\circ\Phi\}  \big(X^{\Psi(\x)}(t)\big)\Big]\\
 &= \bar P_{t+s}^{\a,\b} (f\circ\Phi)\big(\Psi(\x)\big) =P_{t+s}^{\a,\b} f(\x),\ \ s,t>0, f\in \B_b(\bar\DD_\infty), \x\in\bar\DD_\infty.\end{split}\end{equation*} So,  $Y^\x(t)$ is a Markov process on $\bar\DD_\infty$ for any $\x\in \bar\DD_\infty$. Moreover, as shown in  \cite[\S 3]{FW}  that $P_t^{\a,\b}$ is associated to the  symmetric Dirichlet form
$(\EE_{\a,\b},\D(\EE_{\a,\b}))$ on $L^2(\Xi_{\a,\b})$.

It is now the position to  state the main results in the paper. Let
\beg{equation*}\beg{split} &K_{a,b} = \ff {1_{[\ff 1 4,\infty)}(a\land b)} 4 \Big(\ss{(4a-1)(4b-1)}+ 2(a+b)-1\Big),\ \ a,b>0;\\
&\rr(s,t)=\int_{s\land t}^{s\lor t} \ff {\d r} {\ss{r(1-r)}},\ \ s,t\in [0,1].\end{split}\end{equation*}

\beg{thm}\label{T1.1} Assume $ a_i\land b_i \ge \ff 1 4$ for all $i\ge 1.$   Then for any positive $f\in \B_b(\bar\DD_\infty)$ and $p>1$, the following Harnack inequality holds:
\beq\label{H} (P_t^{\a,\b} f)^p(\x)\le (P_tf^p(\y)) \exp\bigg[\ff p{p-1}\sum_{i=1}^\infty \ff{\rr(\psi_i(\x), \psi_i(\y))^2K_{a_i,b_i}}{\exp[2K_{a_i,b_i}t]-1}\bigg],\ \ \x,\y\in \bar\DD_\infty, t>0,\end{equation} where
when $K_{a_i,b_i}=0$ we set $\ff{K_{a_i,b_i}}{\exp[2K_{a_i,b_i}t]-1}=\ff 1 {2t}.$\end{thm}

The following  is a consequence of Theorem \ref{T1.1}.

\beg{cor}\label{C1.2'} Assume $a_i\land b_i\ge \ff 1 4$  for large $i\ge 1$. If
\beq\label{C} \lim_{i\to\infty}\ff{a_i+b_i}{\log i}=\infty,\end{equation} then:
\beg{enumerate}\item[$(1)$] $\Xi_{\a,\b}$ is the unique invariant probability measure of $P_t^{\a,\b}$, and for any $t>0$, $P_t^{\a,\b}$ has a symmetric density $p_t^{\a,\b}(\x,\y)$ with respect to
$\Xi_{\a,\b}$ such that
\beq\label{CB} C^{-1}\e^{-c_0\gg(t)}\le \inf_{\x,\y\in \bar\DD_\infty} p_t^{\a,\b}(\x,\y)\le \sup_{\x,\y\in \bar\DD_\infty} p_t^{\a,\b}(\x,\y)\le C\e^{c_0\gg(t)},\ \ \ t>0\end{equation} holds for some constant $C\ge 1$ and $c_0:=2\rr(0,1),$ where
$$\gg(t):= \sum_{i=1}^\infty \ff{K_{a_i,b_i}}{\exp[K_{a_i,b_i}t]-1}<\infty,\ \ t>0.$$ If $\inf_{i\ge 1}(a_i\land b_i)\ge \ff 1 4$ then $\eqref{CB}$ holds for $C=1.$
 \item[$(2)$] $P_t^{\a,\b}$ is strong Feller with respect to the metric
$${\bf d}(\x,\y):= \bigg(\sum_{i=1}^\infty i^{-2}\rr(\psi_i(\x), \psi_i(\y))^2\Big)^{\ff 1 2},\ \ \ \xi,\eta\in \bar\DD_\infty.$$
\item[$(3)$] Let $\ll=\inf_{i\ge 1}(a_i+b_i)$. Then there exists a constant $c>0$ such that
$$\sup_{\x,\y\in \bar\DD_\infty}|p_t^{\a,\b}(x,y)-1| \le c\,\e^{-\ll t},\ \ t> 1.$$ \end{enumerate}\end{cor}

\paragraph{Remark 1.1}  (1) If $\a,\b$ satisfies \eqref{dustless}, then  $\Xi_{\a,\b}$ is fully supported on the simplex $\DD_\infty$, so that due to Corollary \ref{C1.2'}(1) we have $Y(t)\in \DD_\infty\ \P$-a.s. for any $t>0$ and any starting point $Y(0)\in \bar\DD_\infty.$

(2) It is well known that   the uniform heat kernel upper bound $C\e^{c_0\gg(t)}$ of the heat kernel implies the super log-Sobolev inequality (see \cite[Theorem 5.1.7]{Wbook} or \cite[Theorem 2.2.3]{Davies})
\beq\label{SLS} \Xi_{\a,\b}(f^2\log f^2)\le r\EE_{\a,\b}(f,f) + \log C+ c_0\gg(r),\ \ r>0, f\in\D(\EE_{\a,\b}), \Xi_{\a,\b}(f^2)=1,\end{equation}  as well as the super Poincar\'e inequality (see \cite[Theorem 3.3.15]{Wbook} or \cite[Theorem 4.5]{W00b})
\beq\label{SP0} \Xi_{\a,\b}(f^2)\le  r\EE_{\a,\b}(f,f) + \bb(r) \Xi_{\a,\b}(|f|)^2,\ \ \ r>0, f \in\D(\EE_{\a,\b})\end{equation} for
$$\bb(r):=C \inf_{t>0} \ff r t \exp\Big[c_0\gg(t) +\ff t r -1\Big],\ \ r>0.$$
This strengthens the log-Sobolev inequality derived in \cite{FW}.

(3) Theorem \ref{T1.1}(3) is stronger than the uniform ergodicity (also called strong ergodicity):
$$\sup_{\x\in \bar\DD_\infty} \|P_t^{\a,\b}(\x,\cdot)-\Xi_{\a,\b}\|_{var} \le C \e^{-\ll t},\ \  t\ge 0$$ for some constant $C>0$, where $\|\cdot\|_{var}$ is the total variational  and  $$P_t^{\a,\b}(\x,\d \y):= p_t^{\a,\b}(\x,\y)\Xi_{\a,\b}(\d\y)$$ is the transition probability kernel of the infinite-dimensional diffusion process $Y(t)$.

(4) We also like to mention that by using explicit formula of the heat kernel, the super log-Sobolev inequality has been presented in
\cite[Theorem 4.1]{FSWX} for the infinite-many-neutral-alleles diffusion processes associated to the Poisson-Dirichlet distributions, which are the image of the corresponding GEM distributions of the descending order statistic.

\

To illustrate the above results, we consider below a special case where $a_i+b_i\ge b i$ for some constant $b>0.$ This covers the two-parameter GEM   case where
$a_i= \ff{1-\aa} 2$ and $b_i= \ff{\theta +\aa i}2$ for some constants $\aa\in (0,\ff 1 2]$ and $\theta\ge \ff 1 2-\aa.$

\beg{cor}\label{C1.2} Assume $\inf_{i\ge 1} b_i\ge \ff 1 2$, $a_i\ge \ff 1 4$ for large enough $i\ge 1$, and   $a_i+b_i\ge b i$   for some constant $b>0$ and all $i\ge 1.$  Then there exists a constant $c>0$ such that
\beq\label{T1}\e^{-  c t^{-2}}  \le p_t^{\a,\b}\le \e^{ c t^{-2}} ,\ \ \ t>0,\end{equation} and
\beq\label{T2} \sup_{\x,\y\in\bar\DD_\infty} |p_t^{\a,\b}(\x,\y)-1|\le \e^{ ct^{-2}-\ll t},\ \ t>0, \end{equation} where $\ll:=\inf_{i\ge 1} (a_i+b_i).$ Consequently,  $\eqref{SLS}$  with $\gg(r)= \ff C{r^2}$ and  $\eqref{SP0}$   with $\bb(r)= \exp\big[\ff C{r^{2/3}}\big]$ hold for some constant $C>0$.
\end{cor}

The remainder of the paper is organized as follows. In Section 2 we establish the Harnack inequality and super Poincar\'e inequality  for the Wright-Fisher diffusion processes, which are used  in Section 3 to prove    Theorem \ref{T1.1} and Corollaries \ref{C1.2'}-\ref{C1.2}.

\section{Functional inequalities  for  the Wright-Fisher diffusion processes}

For $a,b>0$, consider the following SDEs on $[0,1]$:
\beq\label{2.1} \d x(t)= \{a-(a+b)x(t)\}\d t + \ss{x(t)(1-x(t))}\,\d B(t),\end{equation} where
$B(t)$ is a one-dimensional Brownian motion. Let $P_t^{a,b}$ be the   Markov semigroup of the solution. Then $P_t^{a,b}$ is symmetric with respect to $\pi_{a,b}$ and, see e.g. \cite[\S 9]{QQ},  has a  density $p_t^{a,b}(x,y)$ with respect to $\pi_{a,b}.$

In this section we investigate the Harnack inequality for $P_t^{a,b}$ and the super Poincar\'e inequality for the associated Dirichlet form
$$\EE_{a,b}(f,f):=\ff 1 2 \int_0^1  x(1-x)f'(x)^2 \d\pi_{a,b},\ \ \ f\in \D(\EE_{a,b}),$$ where $\D(\EE_{a,b})$ is the completion of $C_b^1([0,1])$ under the corresponding $\EE_{a,b}^1$-norm.  These inequalities imply heat kernel estimates and will   be applied in the next section to prove Theorem \ref{T1.1} and Corollaries \ref{C1.2'}-\ref{C1.2}.

We will see in Remark 2.1(2) and the proof of Theorem \ref{T2.2} that the Harnack inequality \eqref{H1} we present below  implies the sharp super Poincar\'e inequality for $a\land b\ge \ff 1 4,$ and the sharp super Poincar\'e inequality for $a\land b\le \ff 1 4 $ will be proved using isoperimetric constants.

\subsection{Harnack inequality and heat kernel estimates}
For any $x\in [0,1]$ and $r>0$, let $B_\rr(x,r)=\{y\in [0,1]: \rr(x,y)<r\}.$

\beg{thm}\label{T2.1} Let $a\land b\ge \ff 1 4.$ Then for any $p>1$ and positive $f\in\B_b([0,1])$, the following  Harnack inequality holds:
\beq\label{H1} (P_t^{a,b} f)^p(x)\le (P_t^{a,b}f^p(y))\exp\bigg[\ff{pK_{a,b} \rr(x,y)^2}{(p-1)(\exp[2K_{a,b}t]-1)}\bigg],\ \
x,y\in [0,1], t>0.\end{equation} Consequently, the heat kernel $p_t^{a,b}$ satisfies
\beq\label{H2} \exp\bigg[-\ff{2K_{a,b}\rr(x,y)^2}{\exp[K_{a,b}t]-1}\bigg]\le p_t^{a,b}(x,y) \le \exp\bigg[\ff{2K_{a,b}\rr(0,1)^2}{\exp[K_{a,b}t]-1}\bigg],\ \ t>0, x,y\in [0,1].\end{equation}\end{thm}

\beg{proof} (a) We first observe that \eqref{H2} follows from \eqref{H1}. Let $p=2$ and $P=P_{\ff t 2}$, \eqref{H1} implies
\beq\label{DE}(P f)^2(x)\le (Pf^2(y))\e^{\Psi(x,y)},\ \ x,y\in [0,1], 0\le  f\in \B_b([0,1]),\end{equation} where $\Psi(x,y):= \ff{2K_{a,b} \rr(x,y)^2}{\exp[K_{a,b}t]-1}.$ So, applying \cite[Theorem 1.4.1(5)]{Wbook2} with $\Phi(r)=r^2$ and using the symmetry of $p(x,y):= p_{\ff t 2}^{a,b}(x,y)$, we obtain
$$p_t^{a,b}(x,y)= \int_{[0,1]} p(x,z)p(y,z)\pi_{a,b}(\d z) \ge \e^{-\Psi(x,y)},$$ which implies the desired lower bound estimate in \eqref{H2}.  Next, by
\cite[Theorem 1.4.1(6)]{Wbook2}, \eqref{DE} implies
$$(Pf)^2(x)\le \ff 1 {\int_0^1 \e^{-\Psi(x,y)}\pi_{a,b}(\d y)},\ \ \pi_{a,b}(f^2)\le 1.$$ Taking     $f(z)=\ff{ p(x,z)}{\ss{p_t(x,x)}}$, we arrive at
$$p_t^{a,b}(x,x)^2 =  (P f(x))^2\le \ff 1 {\int_{[0,1]} \e^{-\Psi(x,y)}\pi_{a,b}(\d y)} \le \exp\bigg[\ff{2K_{a,b}\rr(0,1)^2}{\exp[K_{a,b}t]-1}\bigg],\ \ x\in[0,1].$$ This implies  the desired upper bound estimate in \eqref{H2} since
\beg{equation*}\beg{split} & p_{t}^{a,b}(x,y) = \int_{[0,1]}p(x,z)p(y,z)\pi_{a,b}(\d z)\\
&\le \bigg(\int_{[0,1]}p(x,z)^2\pi_{a,b}(\d z)\bigg)^{\ff 1 2}
\bigg(\int_{[0,1]}p(y,z)^2\pi_{a,b}(\d z)\bigg)^{\ff 1 2}=\ss{p_{t}^{a,b}(x,x)p_t^{a,b}(y,y)}.\end{split}\end{equation*}

(b) next, we prove the Harnack inequality \eqref{H1} using coupling by change of measures. Let $T>0$ and $x,y\in [0,1]$ be fixed. Without loss of
generality, we assume that $y>x.$ Let $x(t)$ solve \eqref{2.1} for $x(0)=x$, and let $y(t)$ solve the following equation on $[0,1]$ with reflection with $y(0)=y$:
\beq\label{2.2}  \beg{split} \d y(t)= &\{a-(a+b)y(t)\}\d t + \ss{y(t)(1-y(t))}\,\d B(t)\\
&-1_{[0,\tau)}(t)\ss{y(t)(1-y(t))}\,\xi(t)\d t,\end{split}\end{equation} where $\tau:=\inf\{t\ge 0: x(t)=y(t)\}$ is the coupling time and
$$\xi(t):= \ff{\rr(x,y) \exp[K_{a,b}t]}{\int_0^T \exp[2K_{a,b}t]\d t},\ \ t\ge 0.$$
Below, we prove the   inequality
\beq\label{A*}
 \d\rr(x(t),y(t)) \le -\big\{K_{a,b}\rr(x(t),y(t)) +\xi(t)\big\}\d t,\ \ \ t\in [0,\tau) \end{equation}
by using It\^o's formula for $\rr(x(t),y(t))$,  see \eqref{a1} below. To avoid the singularity of $\rr(x,y)$ for $x<y$ at $x=0$ and $y=1$, one may
prove \eqref{A*}   in a similar way by   applying It\^o's formula to $\rr_\vv(x(t),y(t)):=\int_{x(t)}^{y(t)} \ff{\d s}{\ss{ (s+\vv)(1+\vv-s)}}$ for $\vv>0$ and finally letting $\vv\to 0$.

Obviously, we have $x(t)<y(t)$ for $t<\tau$, and $x(t)=y(t)$ for $t\ge \tau.$ Consequently, $y(t)>0$ and $x(t)<1$ for $t<\tau$.
Therefore, by It\^o's formula we obtain
\beq\label{a1}\beg{split}  &\d\rr(x(t),y(t))\\
&  =\bigg\{\ff{4a-1-2(2a+2b-1)y(t)}{4\ss{y(t)(1-y(t))}}-\ff{4a-1-2(2a+2b-1)x(t)}{4\ss{x(t)(1-x(t))}}
-\xi(t)\bigg\}\d t \end{split}   \end{equation} for $t\in [0,\tau).$  Since $y(t)>x(t)$ for $t\in [0,\tau),$ we have
\beq\label{a2}\beg{split} &\ff{4a-1-2(2a+2b-1)y(t)}{4\ss{y(t)(1-y(t))}}-\ff{4a-1-2(2a+2b-1)x(t)}{4\ss{x(t)(1-x(t))}}\\
&=\ff 1 4 \int_{x(t)}^{y(t)} \ff{\d}{\d s} \bigg(\ff{4a-1-2(2a+2b-1)s}{\ss{s(1-s)}}\bigg)\d s
 =-\ff 1 8 \int_{x(t)}^{y(t)} \ff{4a-1+4(b-a)s}{\{s(1-s)\}^{\ff 3 2}}\,\d s\\
&\le -c \int_{x(t)}^{y(t)} \ff{\d s}{\ss{s(1-s)}}=-c\rr(x(t),y(t)),\ \ t\in [0,\tau),\end{split}\end{equation} where
\beq\label{a3} c:= \inf_{s\in (0,1)} \ff{4a-1+4(b-a)s}{8s(1-s)} =K_{a,b}.\end{equation}
Since  \eqref{a3} is trivial when $a\land b=\ff 1 4$, we only prove it  for $a\land b>\ff 1 4.$  In this case we have $\ff{4a-1+4(b-a)s}{8s(1-s)}\to\infty$ as $s\to 0$ or $1$, so that the inf is reached in $(0,1).$ It is easy to see that
in $(0,1)$ we have $\ff{\d}{\d s} \big(\ff{4a-1+4(b-a)s}{8s(1-s)}\big)=0$ if and only if
$$4(b-a)s^2+2(4a-1)s-(4a-1)=0,$$ so that the inf is reached at
$$s_0= \ff{\ss{(4a-1)(4b-1)} -(4a-1)}{4(b-a)} = \ff{4a-1}{4a-1+\ss{(4a-1)(4b-1)}}.$$ Thus,
$$c= \ff{4a-1+4(b-a)s_0}{8s_0(1-s_0)}=\ff 1 4 \Big(\ss{(4a-1)(4b-1)} +2a+2b-1\Big)=K_{a,b}.$$
Combining \eqref{a1}, \eqref{a2} and \eqref{a3}, we prove \eqref{A*}.  Consequently,
\beg{equation*}\beg{split} \rr(x(t),y(t))&\le\rr(x,y)\e^{-K_{a,b}t} -\int_0^t \e^{-K_{a,b}(t-s)}\xi(s)\d s\\
&= \ff{\rr(x,y) \e^{-K_{a,b}t} \int_t^T\e^{2K_{a,b}s}\d s}{\int_0^T\e^{2K_{a,b}s}\d s},\ \ t\in [0,\tau).\end{split}\end{equation*}
This implies $\tau\le T$,  so that $x(T)=y(T).$

Now, rewrite \eqref{2.2} as
$$\d y(t)= \{a-(a+b)y(t)\}\d t + \ss{y(t)(1-y(t))}\,\d \tt B(t),$$
where, by Girsanov's theorem,
$$\tt B(t):= B(t) - \int_0^{t\land\tau} \xi(s)\d s,\ \ t\ge 0$$ is a one-dimensional Brownian motion under the probability measure
$\d\Q:=R\d\P$ for
$$R:= \exp\bigg[\int_0^\tau \xi(t)\d B(t)-\ff 1 2\int_0^\tau \xi(t)^2\d t\bigg].$$ So, by the weak uniqueness of the solution to \eqref{2.1}, we have
$$P_T^{a,b}f(y)= \E_\Q f(y(T)) = \E[Rf(y(T))].$$
Combining this with $x(T)=y(T)$ observed above, we obtain
\beg{equation*}\beg{split} (P_T^{a,b}f(y))^p &= \big(\E[Rf(y(T))]\big)^p= \big(\E[Rf(x(T))]\big)^p\\
&\le \big(\E f^p(x(T))\big) \big(\E R^{\ff p{p-1}}\big)^{p-1}= (P_T^{a,b}f^p(x))\big(\E R^{\ff p{p-1}}\big)^{p-1}.\end{split}\end{equation*}
This implies \eqref{H1} since, by the definitions of $R,\xi$ and the fact that $\tau\le T$,
\beg{equation*}\beg{split} \E R^{\ff p{p-1}} &\le \e^{\ff p{2(p-1)^2} \int_0^T \xi(t)^2\d t} \E\exp\bigg[\ff p{p-1} \int_0^\tau \xi(t)\d B(t) -\ff{p^2}{2(p-1)^2}\int_0^\tau \xi(t)^2\d t\bigg]\\
&= \exp\bigg[ \ff p{2(p-1)^2\int_0^T \exp[2K_{a,b}t]\d t}\bigg] =\exp\bigg[\ff{pK_{a,b} \rr(x,y)^2}{(p-1)^2(\exp[2K_{a,b} T]-1)}\bigg].
\end{split}\end{equation*}
\end{proof}

\paragraph{Remark 2.1.} (1) From  the proof we see that the condition $a\land b\ge \ff 1 4$ is more or less essential for the desired explicit Harnack inequality using coupling by change of measures. This condition might be dropped using a localization argument as in \cite{ATW09}, which, however,  will lead to a less explicit  Harnack inequality.

(2) We will see in the proof of Theorem \ref{T2.2} below  that the   Harnack inequality \eqref{H1} also implies the   heat kernel upper bound
\beq\label{HE1} \sup_{x,y\in [0,1]} p_t^{a,b} (x,y) \le \ff{c_{a,b} }{t^{2(a\lor b)}},\ \ \ t\in (0,1]\end{equation} for some constant $c_{a,b}>0$,  which is much better than   \eqref{H2} in short time.  Next,   by repeating the argument in the proof of Lemma 2.3 in \cite{GW},  we see that the Harnack inequality \eqref{H1} implies the following Gaussian type upper bound estimate: for any $\dd>2$ there exists a constant $C(\dd)>0$ such that
\beq\label{H'} p_t^{a,b}(x,y)\le \ff{C(\dd)\exp[-\ff{\rr(x,y)^2}{\dd t}+C(\dd)t]}{\ss{\pi_{a,b}(B_{\ss 2 \rr}(x,\ss{t})) \pi_{a,b}(B_{\ss 2 \rr}(y,\ss{t}))}}, \ \ t>0, x,y\in [0,1],\end{equation} where we have used the fact that $\ss 2\rr$, rather than $\rr$, is the intrinsic distance induced by the diffusion process. Moreover, according to \cite[Theorem 7.2]{GR} which works for the present case by using the transform $x\mapsto \ff x 2 +\ff 1 2$ which maps  $[-1,1]$ therein onto the present $[0,1],$ there exists constants $c_1,c_2,c_1',c_2'>0$ such that
\beq\label{H''}\ff{c_1'\exp[-c_2' \ff{\rr(x,y)^2}{t}]}{\ss{\pi_{a,b}(B_{ \rr}(x,\ff {\ss t}{\ss 2})) \pi_{a,b}(B_{\rr}(y,\ff {\ss t}{\ss 2}))}}\le  p_t^{a,b}(x,y)\le \ff{c_1\exp[-c_2 \ff{\rr(x,y)^2}{t}]}{\ss{\pi_{a,b}(B_{\rr}(x,\ff {\ss t}{\ss 2})) \pi_{a,b}(B_{\rr}(y, \ff {\ss t}{\ss 2}))}}\end{equation} holds for all $t\in (0,1], x,y\in [0,1].$
However, all these estimates can not be extended to infinite-dimensions.

(3) The leading term  of the heat kernel $p_t^{a,b}(x,y)$ has been figured out in the last display in \cite[\S 9]{QQ} as follows
for $b_0=2a$ and $b_1=2b$ (since the reference measure used there is $\d y$ rather than the invariant measure $\pi_{a,b}(\d y)$,  we multiply the factor $y^{1-2a}(1-y)^{1-2b}$):
$$p_t^{a,b}(x,y) \sim   \ff 1 {t^{2a}} \e^{-\ff{x_l+y_l}t} \psi_{2a}\Big(\ff{x_ly_l}{t^2}\Big)
 + \ff 1 {t^{2b}} \e^{-\ff{x_r+y_r}t} \psi_{2b}\Big(\ff{x_ry_r}{t^2}\Big),$$ where $\ss{s_l}:=\sin^{-1}\ss s, \ss{s_r}:= \sin^{-1} \ss{1-s}$ for $s\in [0,1]$,  and
 $$\psi_{b}(x):= \sum_{j=0}^\infty \ff{x^j}{j!\, \GG(j+b)},\ \ \ x,b>0.$$ This suggests
 $$\sup_{x,y\in [0,1]} p_t^{a,b}(x,y)\ge \max_{x\in \{t,1-t\}} p_t^{a,b}(x, x)\ge   \ff{c}{t^{(2a)\lor (2b)}},\ \ t\in (0,1]$$ for some constant $c>0$, so that the above uniform heat kernel estimate \eqref{HE1} implied by the Harnack inequality is sharp for $a\land b\ge \ff 1 4$. See Corollary  \ref{C2.3} below for a sharp uniform heat kernel estimate
 also for $a\land b\le \ff 1 4$ using the super Poincar\'e inequality, which is of order $t^{-(2a)\lor(2b)\lor \ff 1 2}$.

\subsection{Super Poincar\'e inequality and   heat kernel estimates}

According to \cite{W00},   the Dirichlet form $(\EE_{a,b}, \D(\EE_{a,b}))$ is said to satisfy the super Poincar\'e inequality if there exists
a function $\bb: (0,\infty)\to (0,\infty)$ such that
\beq\label{SP} \pi_{a,b}(f^2)\le r\EE_{a,b}(f,f)+ \bb(r)\pi_{a,b}(|f|)^2,\ \ r>0, f\in \D(\EE_{a,b}).\end{equation} As $ \D(\EE_{a,b})$ is the closure of $C_b^1([0,1])$ under the associated Dirichlet norm, one only needs to verify the inequality for $f\in C_b^1([0,1]).$

\beg{thm}\label{T2.2} There exists a constant $c=c(a,b)>0$ such that the following super Poincar\'e inequality
\beq\label{SP2}  \pi_{a,b}(f^2)\le r\EE_{a,b}(f,f)+ \Big(1\lor \ff c {r^{\ff 1 2 \lor (2a)\lor(2b)}}\Big)\pi_{a,b}(|f|)^2,\ \ r>0, f\in \D(\EE_{a,b}).\end{equation}
On the other hand,  the super Poincar\'e inequality $\eqref{SP}$ implies \beq\label{OO} \liminf_{r\to 0} \bb(r)r^{\ff 1 2 \lor (2a)\lor(2b)}>0,\end{equation} so that $\eqref{SP2}$  is sharp for small $r>0.$ \end{thm}

\beg{proof}  (1) The proof of \eqref{SP2} consists of the following four steps.

(1a) It is easy to see that the generator $L_{a,b}:= \ff 1 2 x (1-x)\pp_x^2 +(a-(a+b)x)\pp_x$ has a spectral gap $\ll_1= a+b$ with the first eigenvalue
$u(x):= x-\ff{a}{a+b}.$ Then the Poincar\'e inequality
\beq\label{P} \pi_{a,b} (f^2)\le \ff 1 {a+b} \EE_{a,b}(f,f) +\pi_{a,b}(f)^2,\ \ f\in \D(\EE_{a,b})\end{equation} holds. Thus, for the first assertion it suffices to prove \eqref{SP2} for small $r>0$, say $r\in (0,1].$

 (1b) To prove \eqref{SP2} for $r\in (0,1],$ we first consider $a\land b\ge \ff 1 4$ and prove \eqref{HE1} using the Harnack inequality \eqref{H1}. Since $K_{a,b}\ge 0$, we have
 $$\ff{K_{a,b}}{\exp[2K_{a,b} t]-1}\le \ff 1 {2 t},\ \ \ t>0.$$ So, by    \eqref{H1} with $p=2$ we obtain
$$(P_t^{a,b} f)^2(x) \exp\Big[-\ff{\rr(x,y)^2}t\Big]\le P_t^{a,b} f^2(y),\ \ \ t>0, x,y\in [0,1], f\in \B_b([0,1]).$$ Let $B_\rr(x,r)=\{y\in [0,1]: \rr(x,y)\le r\}$ for $x\in [0,1]$ and $r>0$. This implies
\beg{equation*}\beg{split} (P_t^{a,b} f)^2(x) \e^{-2} \pi_{a,b}\Big(B_\rr\big(x,\ss{2t}\big)\Big)&\le (P_t^{a,b} f)^2(x)\int_{B_\rr(x,\ss{2t})}  \exp\Big[-\ff{\rr(x,y)^2}t\Big]\pi_{a,b}(\d y)\\
& \le \int_0^1 P_t^{a,b}f^2(y)\pi_{a,b}(\d y) \le 1,\ \ \pi_{a,b}(f^2)\le 1.\end{split}\end{equation*} Taking
$$f(z)= \ff{p_t^{a,b}(x,z)}{\ss{p_{2t}^{a,b}(x,x)}},\ \ z\in [0,1],$$   we arrive at
\beq\label{DD0}p_{2t}^{a,b}(x,x)\le \ff{\e^2}{\pi_{a,b}(B_\rr(x,\ss{2t}))},\ \ \ x\in [0,1], t>0.\end{equation}
Similar but less explicit estimates can be derived from \eqref{H'} or \eqref{H''}. We intend to prove
\beq\label{D1} \inf_{x\in [0,1]}\pi_{a,b}\Big(B_\rr\big(x,\ss{t}\big)\Big)\ge c_0 t^{2(a\lor b)},\ \ \ t\in [0,1]\end{equation} for some constant $c_0>0$, so that \eqref{HE1} follows from \eqref{DD0}.

Let $x\in [0,\ff 1 2]$ and take $t_0= \rr(\ff 1 2, \ff 3 4)^2$. Then there exists a unique $y_t\in (x,\ff 3 4]$ such that
$$\ss{t\land t_0}=\rr(x,y_t) =\int_x^{y_t} \ff{\d s}{\ss{s(1-s)} } \le 2 \int_x^{y_t}\ff {\d s}  {\ss s} = 4\big(\ss{y_t}-\ss x\big).$$ So,
$$B_\rr\big(x,\ss t\big)\supset B_\rr\big(x,\ss{t_0\land t}\big)\supset \big[x, \big(\ss x +\ff 1 4 \ss{t\land t_0}\big)^2\big].$$
Combining this with $\big(\ss x +\ff 1 4\ss{t\land t_0}\big)^2 \le y_t^2 \le \ff 3 4,$ and noting that $4a\ge 1$, we obtain
\beg{equation*} \beg{split} & \int_{B_\rr(x,\ss t)} s^{2a-1}(1-s)^{2b-1}\d s \ge \ff 1 {4^{(2b-1)^+}}\int_{(\ss x+ \ff 1 8 \ss{t\land t_0})^2}^{(\ss x+ \ff 1 4 \ss{t\land t_0})^2}s^{2a-1}\d s \\
&=  \ff 1 {2a 4^{(2b-1)^+}}\Big(\big(\ss x +  \ff 1 4 \ss{t\land t_0}\big)^{4a}- \big(\ss x +  \ff 1 8 \ss{t\land t_0}\big)^{4a}\Big)\\
&\ge  \ff {\ss{t\land t_0}} { 4^{1+(2b-1)^+}}\big(\ss x +  \ff 1 8 \ss{t\land t_0}\big)^{4a-1}\ge \ff {(t\land t_0)^{2a}} { 8^{4a-1}4^{1+(2b-1)^+}}.\end{split}\end{equation*} Therefore,
$$\inf_{x\in [0,\ff 1 2]} \pi_{a,b}\Big(B_\rr\big(x,\ss{t}\big)\Big)\ge c_1 t^{2a},\ \ \ t\in [0,1]$$ holds for some constant $c_1>0.$ Similarly, we have
$$\inf_{x\in [\ff 1 2, 1]} \pi_{a,b}\Big(B_\rr\big(x,\ss{t}\big)\Big)\ge c_2 t^{2b},\ \ \ t\in [0,1]$$ for some constant $c_2>0.$ Combining them together we prove \eqref{D1}, and hence \eqref{HE1} as observed above.

Now, according to \cite[Theorem 3.3.15]{Wbook} or \cite[Theorem 4.5]{W00b}, \eqref{HE1} implies the super Poincar\'e inequality \eqref{SP} for
$$\bb(r):= \inf_{t>0} \Big\{\ff r t \e^{\ff t r-1} \|p_t^{a,b}\|_\infty\Big\} \le \|p_r^{a,b}\|_\infty \le \ff{c_0}{ r^{2(a\lor b)}},\ \ r\in (0,1].$$ That is, \eqref{SP2} holds for $r\in (0,1].$

(1c) Next, we consider the case that $a\lor b\le \ff 1 4$, and prove \eqref{SP} for small $r>0$  using  isoperimetric constants. Let $\mu_a(\d x)= 1_{[0,\ff 1 2]}(x) x^{2a-1}\d x. $ Let $\mu_a^\pp$ be the boundary measure induced by $\mu_a$ under the intrinsic metric $\rr$. We have
\beq\label{F1} \mu_a^\pp(\{x\}) :=\lim_{\vv\to 0} \ff{\mu_a(\{y:\rr(y,x)\le \vv\})}{2\vv} \ge c_1 x^{2a-\ff 1 2},\ \ \ x\in (0,  1/2)\end{equation} for some constant $c_1>0.$
Now, for any
set $A\subset [0,\ff 1 2]$ with $\mu_a(A)\in (0, \mu_a([0,\ff 1 2]))$, let $\pp_0 A$ be the set of boundary points of $A$ included in $(0,\ff 1 2).$ Then $\pp_0 A\ne\emptyset.$ It follows from \eqref{F1} and $2a-\ff 1 2\le 0$ that
 $$k(s):= \inf_{A\subset [0,\ff 1 2], \mu_a(A)\le s} \ff{\mu_a^\pp(\pp_0 A)}{\mu_a(A)} \ge \ff {c_1} s,\ \ \ 0< s< \mu_a([0,1/2]).$$ So,
\beq\label{DJ}k^{-1}(2 r^{-1/2})\ge c_2 \ss r\end{equation} holds for some constant $c_2>0$ and small $r>0.$ Therefore, according to \cite[Theorem 3.4.16]{Wbook}, the super Poincar\'e inequality
\beq\label{SP3} \mu_{a}(f^2) \le   r  \int_0^{\ff 1 2} x f'(x)^2\mu_{a}(\d x) + \ff c {\ss r} \mu_{a} (|f| )^2,\ \ r\in (0,1], f\in C_b^1([0,1/2])\end{equation}  holds for some constant $c>0$.  In case the book \cite{Wbook} is not easy to find, we present below a brief proof of the assertion, see also the proof of Theorem 3.4(1) in \cite{W00} where the last term in the first display should be changed into $\ff {2k(r)}r.$  In fact, let $f\in C_b^1([0,\ff 1 2])$ with $\mu_a(|f|)=1$. We have $\mu_a(f^2>t)\le t^{-1/2}$ so that by the coarea formula,
\beg{equation*}\beg{split} &  \int_0^{\ff 1 2} \ss{x(1-x)}\,|(f^2)'(x)| \,\mu_\aa(\d x)  =\int_0^\infty \mu_a^\pp \big(\{f^2=t\}\setminus\{0,1/2 \}\big) \d t\ge k(s)\int_{s^{-2}}^\infty \mu_a(f^2>t)\d t\\
& \le k(s)\mu_a(f^2) -k(s) \int_0^{s^{-2}} \ff{\d t}{\ss t}=k(s)\mu_a(f^2)-\ff{2 k(s)} s,\ \ s\in (0, \mu_a([0,1/2]).\end{split}\end{equation*} Combining this with \beg{equation*}\beg{split} &  \int_0^{\ff 1 2} \ss{x(1-x)}|(f^2)'(x)| \mu_\aa(\d x) \le 2\ss{\mu_a(f^2)} \bigg(\int_0^{\ff 1 2} x (1-x) f'(x)^2 \mu_a(\d x)\bigg)^{\ff 1 2} \\
&\le \ff 2 {k(s)} \int_0^{\ff 1 2} x (1-x) f'(x)^2 \mu_a(\d x)+\ff {k(s)} 2\mu_a(f^2),\end{split}\end{equation*}  we prove
$$\mu_a(f^2)\le \ff 4 {k(s)^2} \int_0^{\ff 1 2} x (1-x) f'(x)^2 \mu_a(\d x)+\ff 4  s,\ \ \mu_a(|f|)=1, s\in (0,\mu_a([0,1/2]).$$ Taking
$s=k^{-1}(2 r^{-1/2})$ in this inequality  and using \eqref{DJ}, we prove \eqref{SP3} for small $r>0.$  Consequently,
 \beq\label{SP4'} \pi_{a,b}(f^21_{[0,\ff 1 2]}) \le r \int_0^{\ff 1 2}  (1- x ) f'(x)^2\pi_{a,b}(\d x) + \ff c {\ss r} \pi_{a,b} (|f| 1_{[0,\ff 1 2]})^2,\ \ r\in (0,1]\end{equation} holds for some constant $c>0$ and all $f\in C_b^1([0,1]).$

 Similarly, when $b\le \ff 1 4$, we have
\beq\label{SP4} \pi_{a,b}(f^21_{[\ff 1 2,1]}) \le r \int_{\ff 1 2}^1 (1- x ) f'(x)^2\pi_{a,b}(\d x) + \ff c {\ss r} \pi_{a,b} (|f| 1_{[\ff 1 2,1]})^2,\ \ r\in (0,1]\end{equation}  for some constant $c>0$
and all $f\in C_b^1([0,1]).$ Combining them together we prove \eqref{SP2} with $r\in (0,1]$ for $a\lor b\le \ff 1 4.$

(1d) Finally, let $a\land b<\ff 1 4$ but $a\lor b \ge \ff 1 4,$ for instance, we assume that $a<\ff 1 4$ and $b\ge \ff 1 4.$  In this case we have $x^{2a-1}\ge x^{-\ff 1 2}$ for $x\in (0,1]$, but
$x^{2a-1}\le 2^{\ff 1 2-2a} x^{-\ff 1 2}$ for $x\in [\ff 1 2,1].$ So,  by \eqref{SP2} for $a\land b\ge \ff 1 4$ we obtain
\beg{equation*}\beg{split} \pi_{a,b}(f^2 1_{[\ff 1 2,1]}) &\le c_1 \pi_{\ff 1 4, b}(f^2) \le r \EE_{\ff 1 4,b} (f,f) +\ff{c_2}{r^{2b}}\pi_{\ff 1 4,b}(|f|)^2\\
&\le c_3\EE_{a,b}(f,f) + \ff{c_3}{r^{2b}}\pi_{a,b}(|f|)^2,\ \ r\in (0,1], f\in C_b^1([0,1]) \end{split}\end{equation*} for some constants $c_1, c_2 ,c_3>0.$ Combining this with \eqref{SP4'}, we prove \eqref{SP2} for $r\in (0,1]$.

(2) To prove the second assertion, let \eqref{SP} hold for some $b$. Take $f(x)= (\vv-x)^+$ for $vv\in (0,\ff 1 2)$. Then there exists constants $c_1,c_2>0$ such that
$$\pi_{a,b}(f^2)\ge c_1 \vv^{2a+2},\ \ \pi_{a,b}(f) + \EE_{a,b}(f,f)\le c_2 \vv^{2a+1},\ \ \vv\in (0,1/2).$$ So, by \eqref{SP} we obtain
$$\bb(r)\ge \ff 1 {c_2^2}\sup_{\vv\in (0,\ff 1 2)} \Big(\ff{c_1}{\vv^{2a}}-\ff{r c_2}{\vv^{2a+1}}\Big) = c_3 r^{-2a} $$ for some constant $c_3>0$ and small $r>0$. Therefore,
$$\liminf_{r\to 0} \bb(r) r^{2a} \ge c_3>0.$$ Similarly, by taking $f(x)= (x+\vv-1)^+$ in \eqref{SP} we obtain $\liminf_{r\to 0} \bb(r) r^{2b} >0;$ while \eqref{SP} with $f(x):= (x+\vv-\ff 1 2)^+\land (\ff 1 2 +\vv -x)^+ $ implies  $\liminf_{r\to 0} \bb(r) r^{\ff 1 2} >0.$ In conclusion, \eqref{OO} holds.
\end{proof}

We would like to indicate that when   $a\land b>\ff 1 4$,   the desired super Poincar\'e inequality  can also be proved using isoperimetric constants. However,
the argument we used is more straightforward and it
stresses   the sharpness of the   Harnack inequality \eqref{H1}.

\beg{cor}\label{C2.3} There exist constants $c_1,c_2>0$ such that
\beq\label{HB1} \sup_{x,y\in [0,1]} |p_t^{a,b}(x,y)-1|\le \ff{c_1\e^{-(a+b)t}}{(t\land 1)^{\ff 1 2\lor (2a)\lor(2b)}} ,\ \ t>0 ;\end{equation}
\beq\label{HB2} \sup_{x,y\in [0,1]} p_t^{a,b}(x,y)\ge \ff{c_2}{t^{\ff 1 2\lor (2a)\lor(2b)}},\ \ \ t\in (0,1].\end{equation} \end{cor}

\beg{proof} (1) Proof of \eqref{H1}. By \cite[Theorem 3.3.15 (2)]{Wbook} or \cite[Theorem 4.5]{W00b}, \eqref{SP2} implies
\beq\label{CCV} \sup_{x,y\in [0,1]} p_t^{a,b} (x,y)\le \ff{c }{t^{\ff 1 2\lor (2a)\lor(2b)}},\ \ \ t\in (0,2]\end{equation}  for some constant $c>0.$ So, it suffices to prove \eqref{HB1} for $t\ge 2.$
By the Poincar\'e inequality \eqref{P}, we have
$$\|P_t^{a,b} -\pi_{a,b}\|_{2\to 2} \le \e^{-(a+b)t},\ \ \ t\ge 0,$$ where, for any $p,q\ge 1,$  $\|\cdot\|_{p\to q}$ stands for the operator norm from $L^p(\pi_{a,b})$ to
$L^q(\pi_{a,b}).$  Combining this with \eqref{CCV} we obtain
$$\|P_t^{a,b} -\pi_{a,b}\|_{2\to\infty} \le \|P_1^{a,b} -\pi_{a,b}\|_{2\to\infty} \|P_{t-1}^{a,b} -\pi_{a,b}\|_{2\to 2} \le C \e^{-(a+b)t},\ \ t\ge 1$$ for some constant $C>0.$
Therefore, by the symmetry of $P_t^{a,b}$ in $L^2(\pi_{a,b})$, this implies
\beg{equation*}\beg{split}  \sup_{x,y\in [0,1]}p_t^{a,b}(x,y)-1| &= \|P_t^{a,b} -\pi_{a,b}\|_{1\to\infty} \le \|P_{\ff t 2}^{a,b} -\pi_{a,b}\|_{1\to 2}\cdot
 \|P_{\ff t 2}^{a,b} -\pi_{a,b}\|_{2\to \infty}\\
 &= \|P_{\ff t 2}^{a,b} -\pi_{a,b}\|_{2\to \infty}^2 \le C^2 \e^{-(a+b)t},\ \ t\ge 2.\end{split}\end{equation*}
Therefore, \eqref{HB1} holds also for $t\ge 2$.

(2) To prove \eqref{HB2}, we use again \cite[Theorem 3.3.15]{Wbook} or \cite[Theorem 4.5]{W00b}  that \eqref{SP} holds for
$$\bb(r) :=\inf_{t>0} \Big\{\ff r t \e^{\ff t r-1} \|p_t^{a,b}\|_\infty\Big\} \le \|p_r^{a,b}\|_\infty.$$ Combining this with the second assertion in Theorem \ref{T2.2}, we obtain
$$\liminf_{t\to 0} \|p_t^{a,b}\|_\infty t^{2(a)\lor(2b)\lor\ff 1 2} >0,$$  which implies \eqref{HB2} for some constant $c_2>0.$
\end{proof}

\section{Proofs of Theorem \ref{T1.1} and Corollaries \ref{C1.2'}-\ref{C1.2}}

\beg{proof}[Proof of Theorem \ref{T1.1}]   By \eqref{F2} and Theorem \ref{T2.1}, we have
\beg{equation}\label{JE}\beg{split}& (P_t^{\a,\b}f(\x))^p =\bar P_t^{\a,\b}(f\circ\Phi)(\Psi(x))= \Big\{\Big(\prod_{i=1}^\infty P_t^{a_i,b_i}\Big)f(\Psi(\x))\Big\}^p\\
&\le \Big\{\Big(\prod_{i=1}^\infty  P_t^{a_i,b_i}\Big)f^p(\Psi(\y))\Big\}\exp\bigg[\sum_{i=1}^\infty \ff{pK_{a_i,b_i}\rr(\psi_i(\x),\psi_i(\y))^2}{(p-1)(\exp[2K_{a_i,b_i}t]-1)}\bigg]\\
&= (P_t^{\a,\b} f^p(\y)) \exp\bigg[\ff p{p-1} \sum_{i=1}^\infty \ff{K_{a_i,b_i}\rr(\psi_i(\x),\psi_i(\y))^2}{\exp[2K_{a_i,b_i}t]-1}\bigg].\end{split}\end{equation} Thus,  \eqref{H} holds.\end{proof}

\beg{proof}[Proof of Corollary \ref{C1.2'}]  (a) Let $i_0\ge 0$ such that $a_i\land b_i\ge \ff 1 4$ for $i> i_0.$ It is easy to see that   \eqref{C} implies $\gg(t)<\infty$ for all $t>0.$ By \eqref{HB1}, there exists a constant $c\ge 1$ such that
$$c^{-1} \e^{-c_0t^{-1}}\le p_t^{a_i,b_i}(x,y)\le c \e^{c_0 t^{-1}},\ \ t>0, 1\le i\le i_0.$$ Combining this with
 \eqref{H2} we obtain
$$C^{-1} \e^{-c_0 \gg(t)}\le \prod_{i=1}^\infty p_t^{a_i,b_i} (x_i,y_i)\le C\e^{c_0\gg(t)},\ \ t>0, \x=(x_1, x_2, \cdots), \y=(y_1,y_2,\cdots)\in [0,1]^\N$$ for some constant $C\ge 1,$ and $C=1$ if $i_0=0.$ So, according to \eqref{F} and the definition of $\Xi_{\a,\b}$, the density of $P_t^{\a,\b}$ with respect to $\Xi_{\a,\b}$ exists and is given by
$$p_t^{\a,\b}(\x,\y)=  \prod_{i=1}^\infty p_t^{a_i,b_i} (\psi_i(\x),\psi_i(\y)),\ \ t>0,\x,\y\in \bar\DD_\infty.$$
Thus, \eqref{CB} holds.

Next, since $P_t^{\a,\b}$ is symmetric in $L^2(\Xi_{\a,\b}),$ $\Xi_{\a,\b}$ is its   invariant probability measure. Moreover, \eqref{CB} implies the Harnack inequality
$$(P_t^{\a,\b}f(\x))^p\le  (P_t^{\a,\b}f(\y))^p C^{2p}\e^{2p c_0 \gg(t)},\ \ \ t>0, p>1, \x,\y\in \bar\DD_\infty$$ for all positive $f\in \B_b(\bar\DD_\infty).$ Then, according to \cite[Theorem 1.4.1(3)]{Wbook2} or \cite[Proposition 3.1(3)]{WY11}, $P_t^{\a,\b}$ has  a unique invariant probability measure. Therefore, the proof of (1) is finished.

(b) Since for any $a,b>0$ the semigroup $P_t^{a,b}$ has a continuous density with respect to $\pi_{a,b}$ (see \cite[\S 9]{QQ} as mentioned in Remark 2.1), it is strong Feller with respect to the metric $\rr$. So, due to the first equality in \eqref{JE}, the strong Feller property of $P_t^{\a,\b}$ with respect to $\bf d$ is not affected by changing finite many $(a_i,b_i)$. Thus, without loss of generality, we may and do assume that $a_i\land b_i\ge \ff 1 4$ for all $i\ge 1.$ In this case,  \eqref{C} implies
$$C(t):= \sup_{i\ge 1}  \ff{i^2K_{a_i,b_i}}{\exp[2K_{a_i,b_i}t]-1}<\infty,\ \ t>0.$$ Then \eqref{H} yields
$$(P_t^{\a,\b}f(\x))^p\le (P_t^{\a,\b} f^p(\y)) \exp\Big[\ff{pC(t) {\bf d} (\x,\y)^2}{p-1}\Big],\ \ t>0, \x,\y\in\bar \DD_\infty.$$
According to \cite[Theorem 1.4.1(1)]{Wbook2} or \cite[Proposition 3.1(1)]{WY11}, $P_t^{\a,\b}$ is strong Feller with respect to the metric ${\bf d},$ i.e.  (2) holds.

(c) Finally, as in  \cite[Theorem 3.1]{FW}, the Poincar\'e inequality
$$\Xi_{\a,\b}(f^2)\le \ff 1 \ll \EE_{\a,\b}(f,f),\ \  f\in \D(\EE_{\a,\b}), \Xi_{\a,\b} (f)=0$$ holds. So,
$$ \|P_t^{\a,\b}-\Xi_{\a,\b}\|_2 \le \e^{-\ll t},\ \ t\ge 0, $$ where $\|\cdot\|_2$ is the $L^2$-norm with respect to $\Xi_{\a,\b}.$ On the other other hand, (1) implies $$\|P_{\ff 1 2}^{\a,\b}-\Xi_{\a,\b}\|_{1\to\infty}  <\infty.$$
Moreover, by the symmetry of $P_{\ff t 2}^{\a,\b}$, we have
$$\|P_{\ff t 2}^{\a,\b}-\Xi_{\a,\b}\|_{L^1(\Xi_{\a,\b})\to L^2(\Xi_{\a,\b})}= \|P_{\ff t 2}^{\a,\b}-\Xi_{\a,\b}\|_{L^2(\Xi_{\a,\b})\to L^\infty(\Xi_{\a,\b})}\le  \|P_{\ff t 2}^{\a,\b}-\Xi_{\a,\b}\|_{2\to \infty},$$ where
$\|\cdot\|_{2\to\infty}$ is defined as $\|\cdot\|_{1\to\infty}$ using $\Xi_{\a,\b}(f^2)\le 1$ in place of $\Xi_{\a,\b}(|f|)\le 1$, i.e. for a linear operator $P$ on $L^2(\Xi_{\a,\b})$,
$$\|P\|_{2\to\infty}:=\sup_{\Xi_{\a,\b}(f^2)\le 1} \sup_{\x\in\bar\DD_\infty} |Pf(\x)|.$$
Therefore,
\beg{equation*}\beg{split} &\|P_t^{\a,\b} -\Xi_{\a,\b}\|_{1\to\infty} \le \|P_{\ff t 2}^{\a,\b}-\Xi_{\a,\b}\|_{L^1(\Xi_{\a,\b})\to L^2(\Xi_{\a,\b})}
\| P_{\ff t 2}^{\a,\b}-\Xi_{\a,\b}\|_{2\to\infty}\\
&\le \| P_{\ff t 2}^{\a,\b}-\Xi_{\a,\b}\|_{2\to\infty}^2\le \| P_{\ff 1 2}^{\a,\b}-\Xi_{\a,\b}\|_{2\to\infty}^2\| P_{\ff {t-1} 2}^{\a,\b}-\Xi_{\a,\b}\|_{2}^2\le c\e^{-\ll t}\end{split}\end{equation*} holds for some constant $c>0$ and all $t\ge 1.$ Now, for any $t\ge 1$ and $\vv>0$, with $f({\mathbf z}):= p_\vv^{\a,\b}({\mathbf z},\y)$ this implies
$$p_{t+\vv}^{\a,\b} (\x,\y)= (P_t^{\a,\b} f)(\x)\le c \e^{-\ll t} \Xi_{\a,\b}(f)= c\e^{-\ll t},\ \ \x,\y\in \bar \DD_\infty.$$ Therefore,
  the proof of (3) is finished.
\end{proof}

\beg{proof}[Proof of Corollary \ref{C1.2}]  By Corollary \ref{C1.2'}(3) and  Remark 1.1(2), it suffices to prove $\gg(t)\le \ff c {t^2}$ for some constant $c>0$ and $t\in (0,1]$ as in this case the   estimate holds for all $t>0$, so that
$$\inf_{t>0} \bigg(\ff r t\exp\Big[c_0\gg(t)+\ff t r-1\Big]\bigg)\le\inf_{t>0} \bigg(\ff r t\exp\Big[\ff{c_0c}{t^2}+\ff t r-1\Big]\bigg)
 \le \exp\Big[\ff C{r^{\ff 2 3}}\Big],\ \ \ r>0$$ holds for some constant $C>0$  by taking   $t=r^{\ff 1 3}$ for $r<1.$

 By the definition of $K_{a_i,b_i}$ and the condition $a_i+b_i\ge b i$ for some constant $b>0$, there exist $i_0\in\N$ and a constant $c_1>0$ such that $K_{a_i,b_i}\ge c_1 i$ for $i\ge i_0.$ Since $K_{a_i,b_i}\ge 0$ and $\e^s\ge \ff s 2 \e^{s/2}$ holds for $s\ge 0$, we get
 \beg{equation*}\beg{split} \gg(t) &\le \ff {i_0} t+ \sum_{i>i_0} \ff{c_1 i}{\exp[c_i i t]-1} \le \ff{i_0} t +\sum_{i>i_0} \ff 2 t \e^{-c_1it/2}\\
 &\le \ff {i_0}t +\ff 2 t \int_1^\infty \e^{-c_1ts/2}\d s =\ff{i_0} t+ \ff 4{c_1t^2}\le \ff c{t^2},\ \ t\in (0,1]\end{split}\end{equation*} for some constant $c>0.$ The proof is finished.
\end{proof}

\beg{thebibliography}{99}

 \bibitem{ATW} M. Arnaudon, A. Thalmaier, F.-Y. Wang,
  \emph{Harnack inequality and heat kernel estimates
  on manifolds with curvature unbounded below,} Bull. Sci. Math.   130(2006), 223--233.

\bibitem{ATW09} M. Arnaudon,  A. Thalmaier,   F.-Y. Wang,
 \emph{Gradient estimates and Harnack inequalities on non-compact Riemannian manifolds,}  Stoch. Proc. Appl. 119(2009), 3653-3670.

\bibitem{ConMoi69} R. J. Connor,  J. E. Moismann,
 \emph{Concepts of independence for proportions with a generalization of the Dirichlet distribution,}  J. Amer. Statist. Assoc. 64(1969), 194-206.

 \bibitem{Davies} E. B. Davies, \emph{Heat Kernels and Spectral Theory,}
Cambridge: Cambridge Univ. Press, 1989.

\bibitem{QQ} C. L. Epstein, R. Mazzeo\emph{ Wright-Fisher diffusion in one dimension,}  SIAM J. Math. Anal.   42(2010),  568--608.

\bibitem{Ewens04} W.J. Ewens, \emph{Mathematical Population Genetics, Vol I} Springer-Verlag, New York, 2004.

 \bibitem{FW}  S. Feng, F.-Y. Wang,  \emph{A class of infinite-dimensional diffusion processes with connection to population genetics,} J. Appl. Probab. 44(2007), 938--949.

 \bibitem{FSWX}  S. Feng, W. Sun, F.-Y. Wang, F. Xu, \emph{Functional inequalities for the two-parameter extension of the infinite-many-neutral-alles diffusion,} J. Funct. Anal. 260(2011), 399--413.

 \bibitem{GW} F.-Z. Gong, F.-Y. Wang, \emph{Heat kernel estimates with application to compactness of manifolds,} Quart. J. Math. 52(2001), 171--180.

 \bibitem{GR} T. Coulhon, G. Kerkyacharian, P. Petrushev, \emph{Heat kernel generated frames in the setting
of Dirichlet spaces,} J Fourier Anal. Appl. 18(2012), 995--1066.

 \bibitem{IW} N. Ikeda, S. Watanabe, \emph{Stochastic Differential Equations and Diffusion Processes,} 2nd Ed. North-Holland, Amsterdam, 1989.

 \bibitem{IshJames01} H. Ishwaran,  L. F. James,
 \emph{Gibbs sampling methods for stick-breaking priors,}  J. Amer. Statist. Assoc. 96(2001), 161-173.

\bibitem{W97} F.-Y. Wang, \emph{Logarithmic Sobolev inequalities on noncompact Riemannian manifolds,} Probab. Theory Relat. Fields
109(1997),  417--424.

\bibitem{W00} 	F.-Y. Wang, \emph{ Functional inequalities for empty essential spectrum,}  J. Funct. Anal. 170(2000), 219--245.

\bibitem{W00b} F.-Y. Wang, \emph{Functional inequalities, semigroup properties and spectrum estimates,}  Infinite Dimensional Analysis, Quantum Probability and Related Topics 3(2000), 263--295.

\bibitem{WY11}  F.-Y. Wang, C. Yuan, \emph{ Harnack inequalities for functional SDEs with multiplicative noise and applications,} Stoch. Proc. Appl. 121(2011), 2692--2710.

\bibitem{Wbook} F.-Y. Wang, \emph{Functional Inequality, Markov Semigroups, and Spectral Theory,} Science Press, 2005, Beijing.

\bibitem{Wbook2} F.-Y. Wang, \emph{Harnack Inequalities and Applications for Stochastic Partial Differential Equations,} Springer, 2013, Berlin.

\end{thebibliography}
\end{document}